\documentclass[12pt]{amsproc}
\usepackage[cp1251]{inputenc}
\usepackage[russian]{babel}
\usepackage{amssymb}

\theoremstyle{plain}

\newtheorem{savkoz}{Савкоз}[section]

\newtheorem{Lemma}[savkoz]{Лемма}
\newtheorem{Statement}[savkoz]{Предложение}
\newtheorem{Theorem}[savkoz]{Теорема}

\newtheorem{Note}[savkoz]{Замечание}

\newcommand{\intl}{\mathop{\int}\limits}

\newcommand{\suml}{\mathop{\sum}\limits}

\renewcommand{\Im}{\operatorname{Im}}

\def\Dom{{\frak D}}
\def\sh{\phantom{.}}
\def\al{\alpha}
\def\la{\lambda}

\def\ch{\operatorname{ch}}
\def\sh{\operatorname{sh}}
\def\eps{\varepsilon}

\def\C{{\mathbb C}}

\newlength{\lenun}
\newlength{\lendu}
\def\Wo{\makebox{\parbox[c][\lendu][b]{\lenun}{\(\stackrel{o}{W}\)
}}}

\makeatletter \@addtoreset{equation}{section} \makeatother

\begin{document}
\title[]{О собственных значениях оператора
Штурма--Лиувилля с потенциалами из пространств Соболева}
\thanks{Работа поддержана грантом РФФИ \No~04-01-00712, грантом
поддержки ведущих научных школ НШ-5247.2006.1 и грантом ИНТАС
\No~05-1000008-7883}
\author[]{Савчук А.~М., Шкаликов А.~А.}
\address{Московский государственный университет
им.~М.~В.~Ломоносова, механико-математический факультет}
\email{artem\_savchuk@mail.ru}
\address{Московский государственный университет
им.~М.~В.~Ломоносова, механико-математический факультет}
\email{ashkalikov@yahoo.com}
\begin{abstract}
В статье изучается асимптотическое поведение собственных значений
оператора Штурма-Лиувилля $Ly= -y'' +q(x)y $ с потенциалами из
соболевского пространства $W_2^{\theta -1}, \theta \geqslant 0$,
включая неклассический случай $\theta \in [0,1)$,  когда потенциал
является распределением.  Результаты получены в новых терминах.
Положим
$$
s_{2k}(q)= \lambda_{k}^{1/2}(q)-k, \quad s_{2k-1}(q)=
\mu_{k}^{1/2}(q)-k-1/2,
$$
где $\{\lambda_k\}_1^{\infty}$ и  $\{\mu_k\}_1^{\infty}$ "---
последовательности собственных значений оператора $L$,
порожденного краевыми условиями Дирихле и Дирихле-Неймана,
соответственно. Построены специальные гильбертовы пространства
$\hat l_2^{\,\theta}$, такие, что отображение $F: W^{\theta-1}_2
\to \hat l_2^{\,\theta}$, определяемое  равенством
$F(q)=\{s_n\}_1^{\infty}$, корректно определено  для всех
$\theta\geqslant 0$. Основной результат заключается в следующем:
при $\theta>0$ отображение $F$ является слабо нелинейным, т.е.
представимо в виде $F(q) =Uq+\Phi(q)$, где $U$ "--- изоморфизм
пространств $W^{\theta-1}_2 $ и $\hat l_2^{\,\theta}$, а $\Phi(q)$
"--- компактное отображение. Более того, доказана оценка
$\|\Phi(q)\|_{\tau} \leqslant C\|q\|_{q-1},$ где значение $\tau
=\tau(\theta)>\theta$ точно указано, а постоянная  $C$  зависит
только от радиуса шара $\|q\|_{\theta} \leqslant R$,  но не
зависит от функции $q$, меняющейся в этом шаре.
\end{abstract}
\begin{flushleft}
УДК~517.984
\end{flushleft}
\maketitle
\bigskip

\section*{\bf \S 1. Введение. Формулировка основной
теоремы.}\refstepcounter{section}

\medskip
Цель настоящей работы -- получить точные асимптотические формулы с
{\it равномерной} оценкой остаточных членов для собственных
значений оператора Штурма--Лиувилля на конечном отрезке. Конечно,
асимптотические формулы для собственных значений этого оператора
хорошо известны с давних времен (см, например, \cite{Ma},
\cite{Le}), но в той форме и в той общности, которая нам нужна,
они не приводились. Мы будем квалифицировать поведение собственных
значений в терминах отображений гильбертовых пространств. Такая
постановка позволит использовать преимущества абстрактной техники.
В частности, точная информация об асимптотике собственных значений
будет получена не только для потенциалов $q$, лежащих в
соболевских пространствах $W_2^{\theta-1}$ с целыми положительными
$\theta$, но и при всех $\theta\geq0$, включая случай
$\theta\in[0,1)$, когда $q$ является распределением. В отличие от
предшествующих работ мы получим равномерные оценки остаточных
членов, когда $q$ меняется в шаре радиуса $R$ пространства
$W_2^{\theta-1}$. Для нелинейных отображений нет аналога теоремы
Банаха--Штейнгауза и получение равномерных оценок требует более
тонкой работы. Имея в виду дальнейшие приложения результатов этой
статьи к {\it обратным задачам Штурма--Лиувилля}, мы ограничимся
рассмотрением операторов, определяемых равенством
\begin{equation} \label{1}
 Ly=-y''+q(x)y,\qquad x\in[0,\pi],
\end{equation}
и краевыми условиями Дирихле, а также краевыми условиями
Дирихле--Неймана, хотя развитые методы позволяют получить
аналогичные результаты для более общих краевых условий.

Напомним определение оператора Штурма--Лиувилля с
потенциалом--раc\-пре\-делением. Пусть $q(x)\in
W_2^{\theta-1}[0,\pi]$, $\theta\in[0,1)$. Тогда существует
первообразная $\sigma(x)=\int q(x)dx\in W_2^{\theta}[0,\pi]$, т.е.
функция $\sigma$, для которой $q(x)=\sigma'(x)$ в смысле
распределений. Известно (см., например \cite[гл. 1]{GS}), что
такая функция $\sigma$ определена однозначно с точностью до
постоянной.

Оператор $L$, порожденный выражением \eqref{1} и краевыми
условиями Дирихле (будем обозначать его через $L_D$) теперь можно
определить следуя нашей работе \cite{SS1} (см. также \cite{SS2},
где приведены другие определения):
\begin{equation} \label{2}
L_Dy=Ly=-(y^{[1]})'-\sigma(x)y^{[1]}-\sigma^2(x)y,\qquad\text{где}\quad
y^{[1]}(x):=y'(x)-\sigma(x)y(x)
\end{equation}
на области

$$
\Dom(L_D)=\{y,\,y^{[1]}\in W_1^1[0,\pi]\vert\ y(0)=y(\pi)=0\}.
$$
Аналогично определяется оператор Дирихле--Неймана: $L_{DN}y=Ly$ на
области
$$
\Dom(L_{DN})=\{y,\,y^{[1]}\in W_1^1[0,\pi]\vert\
y(0)=y^{[1]}(\pi)=0\}.
$$
Конечно, для гладкой функции $\sigma$ правые части \eqref{1} и
\eqref{2} совпадают и мы получаем классические операторы
Штурма--Лиувилля с краевыми условиями Дирихле и Дирихле--Неймана
(во втором случае константу при выборе функции $\sigma$ нужно
взять такой, чтобы $\int_0^\pi q(x)dx=\sigma(\pi)=0$).

Для формулировки основных результатов работы  важную роль играют
пространства $\hat{l}_2^{\,\theta}$, введенные ранее в нашей
работе \cite{SS3}.\footnote{В работе \cite{SS3} для этого
пространства вместо $\hat{l}_2^{\,\theta}$ использовалось
обозначение $l_2^{\ \theta}$, которое трудно отличимо от
обозначения обычного весового пространства $l_2^\theta$.}
Фиксируем $\theta\geqslant0$ и обозначим через $l_2^{\,\theta}$
пространство, состоящее из последовательностей комплексных чисел
$x=\{x_1,x_2,\dots\}$ таких, что
$$
\suml_{k=1}^\infty|x_k|^2k^{2\theta}<\infty.
$$
Введем специальные последовательности
$$
e^{2s-1}=\{k^{-(2s-1)}\}_{k=1}^\infty,\qquad
e^{2s}=\{(-1)^kk^{-(2s-1)}\}_{k=1}^\infty,\quad s=1,2,\dots.
$$
При заданном $\theta\geqslant0$ найдется единственное натуральное
число $m$, такое, что $2m-3/2\leqslant\theta<2m+1/2$. Для этого
числа $\theta $ определим пространство $\hat{l}_2^{\,\theta}$ как
конечномерное расширение пространства $l_2^{\,\theta}$, а именно,
$$
\hat{l}_2^{\,\theta}=l_2^{\,\theta}\oplus span\{e^k\}_1^{2m}.
$$
Тем самым, $\hat{l}_2^{\,\theta}$ состоит из элементов
$\hat{x}=x+\sum_{k=1}^{2m}\al_ke^k$, где $x\in{l}_2^{\,\theta}$ и
$\{\al_k\}_1^{2m}$ -- комплексные числа. Скалярное произведение
элементов $\hat{x},\ \hat{y}\in\hat{l}_2^{\,\theta}$ определяется
формулой
$$
(\hat{x},\hat{y})_\theta=(x,y)_\theta+\suml_{k=1}^{2m}\al_k\overline{\beta_k},
$$
где $(x,y)_\theta$ -- скалярное произведение в ${l}_2^{\,\theta}$,
а $\beta_k$ -- коэффициенты элемента $\hat{y}$ при элементах
$e^k$. В частности $\hat{l}_2^{\,\theta}={l}_2^{\,\theta}$ при
$0\leqslant\theta<1/2$ и
$\hat{l}_2^{\,\theta}={l}_2^{\,\theta}\oplus span\{e^k\}_1^2$ при
$1/2\leqslant\theta<5/2$. Мы будем использовать следующее
утверждение, доказанное в \cite{SS3}.

{\bf Предложение.} {\it При любом $\theta\geqslant0$ оператор
\begin{equation}
\label{T} T\sigma=\{b_k\}_1^\infty,\quad
b_k=\frac2\pi\intl_0^\pi\sigma(x)\sin kxdx
\end{equation}
реализует изоморфизм соболевского пространства $W_2^\theta[0,\pi]$
на $\hat{l}_2^{\,\theta}$.}

Теперь сформулируем основной результат статьи.

{\bf Основная Теорема.} {\it Пусть $\{\la_k\}_1^\infty$ и
$\{\mu_k\}_1^\infty$ -- последовательности собственных значений
операторов $L_D$ и $L_{DN}$, соответственно, с комплексным
потенциалом $q(x)\in W_2^{\theta-1}$ (или $\sigma(x)\in
W_2^\theta$), $\theta\geqslant0$. Образуем последовательность
$\{s_k\}_1^\infty$, где
$$
s_{2j-1}=\sqrt{\mu_j}-(j-1/2),\qquad s_{2j}=\sqrt{\la_j}-j,\quad
j=1,2,\dots.
$$
Здесь значения корня выбираются так, чтобы аргументы  чисел
$\sqrt{\mu_j}$ и $\sqrt{\la_j}$ лежали в промежутке
$(-\pi/2,\,\pi/2]$. Тогда $F(\sigma)=\{s_k\}_1^\infty$ является
отображением из $W_2^\theta$ в $\hat{l}_2^{\,\theta}$, причем
$$
F(\sigma)=-(1/2)T\sigma+\Phi(\sigma),
$$
где линейный оператор $T$ определен \eqref{T}, а нелинейное
отображение $\Phi$ отображает $W_2^\theta$ в $\hat{l}_2^{\,\tau}$,
где
$$
\tau=\left\{\begin{array}{ll} 2\theta,\ \text{если}\
0\leqslant\theta\leqslant1,\\ \theta+1,\ \text{если}\
1\leqslant\theta<\infty.\end{array}\right.
$$
При этом отображение $\Phi:\,W_2^\theta\to\hat{l}_2^{\,\tau}$
является ограниченным в каждом шаре, т.е.
\begin{equation}\label{1.4}
\|\Phi(\sigma)\|_\tau\leq C\|\sigma\|_\theta\quad\text{при всех}\
\sigma,\ \|\sigma\|_\theta\leq R,
\end{equation}
где постоянная $C$ зависит от $R$ но не зависит от $\sigma$ в шаре
$\|\sigma\|_\theta\leq R$.}

Кроме того, в работе будет доказано, что $F(\sigma)$ является
вещественным аналитическим отображением и его производная в точке
$\sigma=0$ совпадает с оператором  $-(1/2)\, T$.

В классическом случае $\theta \geq 0$  наиболее точные асимптотики
для собственных значений при целых  $\theta$  получил Марченко
\cite[гл. 1]{Ma}. Его результат в наших терминах можно
переформулировать так: при $\theta =1,2,\dots,$ отображение
$\Phi(\sigma):W^{\theta}_2 \to \hat{l}_2^{\,\theta +1}$ корректно
определено на каждом элементе. Здесь мы развиваем другую технику,
пригодную для доказательства равномерной ограниченности. Прямая
задача Штурма-Лиувилля с потенциалами-распределениями (случай
$\theta\in [0,1))$ помимо указанных наших статей изучалась
в работах Гринива и Микитюка \cite{HM01}-\cite{HM4}, Каппелера
и Мёра \cite{KM}, Коротяева~\cite{Ko}.

Одним из наиболее важных следствий основной теоремы является
следующий результат: {\it При $\theta>0$ отображение
$\Phi(\sigma):\,W_2^\theta\to\hat{l}_2^{\,\theta}$ является
компактным, т.е. $F(\sigma):\,W_2^\theta\to\hat{l}_2^{\,\theta}$
является слабо нелинейным отображением.} Этот факт будет играть
ключевую роль для получения результатов о глобальной устойчивости
в обратной задаче Штурма--Лиувилля, что мы планируем провести в
другой работе.

В этой статье мы не изучаем случай $\theta=0$. Этот случай особый
и требует отдельного рассмотрения.

\medskip

\section*{\bf \S 2. Случай $0<\theta<1/4$.}\refstepcounter{section}

\medskip

При \(0<\theta<1/4\) основная теорема эквивалентна следующему утверждению:
\begin{Theorem}
\label{thm3.1} Пусть \(\sigma(x)\in W_2^{\theta}\),
$\{\la_n\}_1^\infty$ и $\{\mu_n\}_1^\infty$-- собственные значения
операторов  $L_D$ и $L_{DN}$, соответственно. Тогда
$$
\sqrt{\la_n}=n - (1/2) b_{2n}+\al_{2n},\qquad \sqrt{\mu_n}=n-1/2 -
(1/2) b_{2n-1}+\al_{2n-1},
$$
где
$$
b_{k}=\frac2{\pi}\int\nolimits_0^\pi\sigma(x)\sin(kx)dx,\qquad
\sum_{n=1}^\infty k^{2\theta}|\al_{k}|^2<C,
$$
а постоянная $C$ зависит только от радиуса $R$ шара
$\|\sigma\|_\theta\leqslant R$, в котором находится функция
$\sigma=\int q(x)dx\in L_2[0,\pi]$.
\end{Theorem}
\begin{proof}
Мы будем использовать следующие леммы.

\begin{Lemma}\label{lem3.2}
Пусть число \(\rho\in\mathbb C\) такое, что существует решение
\(\theta=\theta(x,\rho)\) уравнения
\begin{equation}\label{ekv:main}
    \theta(x,\rho)=\rho x+\int\limits_0^x\sigma(t)\,
    \sin(2\theta(t,\rho))\,dt+\dfrac{1}{2\rho}\int\limits_0^x
    \sigma^2(t)\,dt-\dfrac{1}{2\rho}\int\limits_0^x\sigma^2(t)\,
    \cos(2\theta(t,\rho))\,dt.
\end{equation}
Тогда решение \(s(x,\rho)\) уравнения
$$
Ly=\rho^2y,
$$
удовлетворяющее условиям \(s(0,\rho)=0\), \(s^{[1]}(0,\rho)=1\),
допускает представление
\begin{equation}\label{polar}
    \rho s(x,\rho)=r(x,\rho)\,\sin\theta(x,\rho),\\
\end{equation}
где
 $$   r(x,\rho)=\exp\left\{-\int\limits_0^x\sigma(t)\,
    \cos(2\theta(t,\rho))\,dt-\dfrac{1}{2\rho}\int\limits_0^x\sigma^2(t)\,
    \sin(2\theta(t,\rho))\,dt\right\}.
$$
\end{Lemma}
\textit{Доказательство} этого утверждения имеется в п.~2.1
работы~\cite{SS2}. Функции $r$ и $\theta$ названы
модифицированными функциями Прюффера, а
представление~\eqref{polar} --- полярным представлением.
\begin{Lemma}\label{lem3.3}
Пусть \(\sigma(x)\in L_2\), \(\|\sigma\|_{L_2}\leqslant R\), а \
$\nu$ "--- фиксированное число. Положим
$$
\upsilon(x,\rho)=\intl_0^x\sigma(t)\sin(2\rho t)dt+
2\intl_0^x\intl_0^t\sigma(t)\sigma(s)\cos(2\rho t)\sin(2\rho
s)dsdt+\frac1{2\rho}\intl_0^x\sigma^2(t)(1-\cos(2\rho t))dt,
$$
\[
    \Upsilon(\rho)=\max\limits_{0\leqslant x\leqslant\pi}\left|
    \int\limits_0^x\sigma(t)\sin(2\rho t)\,dt\right|+
    \left|\int\limits_0^x\sigma(t)\,\cos(2\rho t)\,dt\right|+
    \dfrac{R^2(1+\varkappa+R\varkappa^2)}{2|\rho|},
\]
где \(\varkappa=\ch 2\pi\nu\). Тогда найдётся абсолютная
постоянная \(\varepsilon>0\) (можно взять \(\varepsilon=2^{-7}\)),
такая, что для всех \(\rho\), лежащих в полосе
\(|\Im\rho|\leqslant\nu\) и удовлетворяющих условию
\begin{align}\label{ekv:ineq}
    \Upsilon(\rho)&<\varepsilon(1+64R^2\varkappa^2)^{-2},
\end{align}
уравнение~\eqref{ekv:main} имеет единственное решение \(\theta(x,\rho)\),
которое представимо в виде
\begin{align}\label{ekv:pred}
    \theta(x,\rho)&=\rho x+f(x,\rho),&
    |f(x,\rho)|&\leqslant C\Upsilon(\rho),
\end{align}
где \(C\) "--- абсолютная постоянная. Более того,
\begin{equation}\label{eq:pred1}
\theta(x,\rho)=\rho x+\upsilon(x,\rho)+r(x,\rho),\qquad
|r(x,\rho)|\leq C(R)\Upsilon^2(\rho).
\end{equation}
\end{Lemma}
\begin{proof}
Это утверждение доказано в лемме~2.1 работы~\cite{SS2}, но при
более ограничительном условии: при \(|\rho|>r_0\), где \(r_0\)
зависит от \(\nu\) и \(\sigma\). Чтобы получить более точную
зависимость только от \(R\) и \(\nu\), проанализируем этапы~1 и~3
доказательства леммы~2.1 из~\cite{SS2}. Здесь мы сокращаем
вычисления, которые легко восстановить самостоятельно или
из~\cite{SS2}.

Обозначим через \(F(\theta)\) правую часть уравнения~\eqref{ekv:main}
и рассмотрим в пространстве \(C[0,\pi]\) отображение
\begin{align*}
    \Phi(\theta)&=F(\theta+\theta_0)-\theta_0,&
    \theta_0&=\theta_0(x,\rho)=\rho x.
\end{align*}
Уравнение~\eqref{ekv:main} перепишем в виде
\begin{align*}
    f&=\Phi(f),&f&=\theta-\theta_0.
\end{align*}
Воспользовавшись тригонометрическими формулами для вычисления
выражения \(F(f+\theta_0)\), получим
\begin{gather*}
    \Phi(f)=\Phi_0(f)+\Phi_1(f)+\Phi_2(f)+\Phi_3(f),\\
    \intertext{где}
    \Phi_0=\int\limits_0^x \sigma(t)\,\sin (2\rho t)\,dt+
    \frac1{2\rho}\intl_0^x\sigma^2(t)(1-\cos(2\rho t))\,dt,    \quad
    \Phi_1(f)=2\int\limits_0^x f(t)\sigma(t)\,\cos (2\rho t)\,dt,\\
    \Phi_2(f)=\int\limits_0^x \sigma(t)\,\cos (2\rho t)\,
    (\sin 2f(t)-2f(t))\,dt-\int\limits_0^x \sigma(t)\,
    \sin (2\rho t)\,(1-\cos 2f(t))\,dt,\\
    \Phi_3(f)=\dfrac{1}{2\rho}\left[\int\limits_0^x
    (1-\cos 2f(t))\cos (2\rho t) \sigma^2(t)\,\,dt+
    \int\limits_0^x\sin 2f(t) \,\sin (2\rho t)\sigma^2(t)\,dt
    \right].
\end{gather*}
Покажем, что отображение \(\Phi^2\) является сжимающим в шаре
достаточно малого радиуса пространства \(C[0,\pi]\). Далее норму в
\(C[0,\pi]\) мы обозначаем так же, как абсолютную величину
\(|\cdot|\). Заметим, что \(\Phi_0\) не зависит от \(f\) и
\begin{equation}\label{ekv:3.0}
    |\Phi_0|<\Upsilon(\rho).
\end{equation}
Отображение \(\Phi_1\) есть линейный оператор, причём
\begin{gather}\notag
    |\Phi_1^2(f)|=4\,\int\limits_0^x f(s)\sigma(s)\,\cos(2\rho
    s)\ \int\limits_s^x\sigma(t)\,\cos(2\rho t)\,dtds.\\
    \intertext{Тогда}\label{ekv:3.1}
    \begin{aligned}|\Phi_1(f)|&<2\sqrt\pi R\varkappa\,|f|,&
    |\Phi_1^2(f)|&<8\sqrt\pi R\varkappa\Upsilon(\rho)|f|\end{aligned}
\end{gather}
(здесь мы учли, что \(\left|\int_s^x\right|<2\Upsilon(\rho)\)). Далее, если
\(f\) и \(g\) лежат в шаре радиуса \(r\) пространства \(C[0,\pi]\), то
\begin{gather}\label{ekv:3.2}
    |\Phi_2(f)-\Phi_2(g)|<4\sqrt\pi R\varkappa\ch 2r\,|f-g|^2<
    8\sqrt\pi R\varkappa\, r\ch2r\,|f-g|,\\
    \label{ekv:3.3}
    |\Phi_3(f)-\Phi_3(g)|\leqslant |\rho|^{-1}R^2\varkappa
    (\sh 2r+\ch2 r)\,|f-g|.
\end{gather}
Последние два неравенства мы получили с учетом следующего
известного факта: \emph{если \(G(\xi)\) "--- аналитическая
функция, то \(|G(\xi)-G(\zeta)|\leqslant M\,|\xi-\zeta|\), где
\(M=\max|G'(\eta)|\) и максимум берётся по \(\eta\), лежащем на
отрезке \([\xi,\zeta]\) комплексной плоскости}.

Положим \(r=1/16\,(1+64R^2\varkappa^2)^{-1}\). Тогда коэффициент
сжатия отображения \(\Phi_1\Phi_2\) не превосходит \( 2^{-5}\).
Из~\eqref{ekv:3.0}--\eqref{ekv:3.3} следует, что если каждое из
чисел \(16\Upsilon R\varkappa \), \(2R^2\varkappa|\rho|^{-1}\),
$8R^3\varkappa^2|\rho|^{-1}$ не превосходит \(2^{-5}\) (а это так,
если выполнено неравенство~\eqref{ekv:ineq} с числом
\(\varepsilon=2^{-7}\)), то каждое из отображений
\(\Phi_j\Phi_k\), \(j,k=0,1,2,3\), является сжимающим с
коэффициентом сжатия \(\leqslant 2^{-5}\). Следовательно,
\(\Phi^2\) "--- сжимающее отображение с коэффициентом сжатия
\(\leqslant 1/2\).

Покажем, что \(\Phi^2\) переводит шар \(|f|\leqslant r\) в себя.
Из определения следует, что \(|\Phi(0)|<\Upsilon=\Upsilon(\rho)\).
Прямыми вычислениями (см.~этап~1 в доказательстве леммы~2.1
работы~\cite{SS2}) можно получить оценку
\begin{equation}\label{f_2}
    |\Phi^2(0)-\upsilon(x,\rho)|< \left[(3+4\sqrt\pi
    R+2R^2)\varkappa^2+2\right]\Upsilon^2(\rho)<\Upsilon(\rho).
\end{equation}
Взяв число \(\varepsilon\) достаточно малым (например,
\(\varepsilon= 2^{-7}\)), с учетом~\eqref{ekv:ineq} получим
\[
    |\Phi^2(0)|<|\Phi^2(0)-\upsilon(x,\rho)|+|\upsilon(x,\rho)|<
    \Upsilon(\rho)+(1+\varkappa R)\Upsilon(\rho)<
    \varepsilon(2+\varkappa R)(1+64R^2\varkappa^2)^{-2}<(1/2)r.
\]
Из неравенства \(|\Phi^2(f)-\Phi^2(0)|<(1/2)|f|\) теперь получаем
\(|\Phi^2(f)|<(1/2)(|f|+r)<r\), если \(|f|\leqslant r\).

Таким образом, \(\Phi^2\) "--- сжимающее отображение в малом шаре пространства
\(C[0,\pi]\). Взяв \(f_0=0\), получим, что последовательные приближения
\(f_k=\Phi(f_{k-1})\) сходятся к решению \(f\) уравнения \(f=\Phi(f)\).
Функция \(\theta=f+\theta_0\) будет решением уравнения~\eqref{ekv:main}.

Функция \(f\) есть сумма ряда
\[
    \Phi(f_0)+[\Phi(f_1)-\Phi(f_0)]+[\Phi^2(f_1)-\Phi^2(f_0)]+\ldots.
\]
Так как \(\Phi^2\) "--- сжатие с коэффициентом \(\leqslant 1/2\), то имеем
оценку
\[
    |f|<|\Phi(f_1)|+\sum\nolimits_{k=1}^{\infty} 2^{-k}|f_1-f_0|
    +\sum\nolimits_{k=1}^{\infty} 2^{-k}|f_2-f_1|<C\Upsilon(\rho),
\]
где \(C\) "--- абсолютная постоянная. Это влечет \eqref{ekv:pred}.

Остается доказать неравенство \eqref{eq:pred1}, эквивалентное
\begin{equation}\label{lt1}
|f-\upsilon|<C(R)\Upsilon^2(\rho).
\end{equation}
Учитывая неравенство
$|f_4-f_3|=|\Phi(f_3)-\Phi(f_2)|<C(R)|f_3-f_2|$ и оценку
\eqref{f_2}, получаем
\begin{multline*}
|f-\upsilon|\leq|f-f_2|+|f_2-\upsilon|<\sum\nolimits_{k=1}^\infty2^{-k}|f_3-f_2|+
\sum\nolimits_{k=1}^\infty2^{-k}|f_4-f_3|+|f_2-\upsilon|<\\
<(1+C(R))|f_3-f_2|+C(R)\Upsilon^2(\rho).
\end{multline*}
Далее, вновь учитывая \eqref{f_2}, имеем
$$
|f_3-f_2|\leq|f_3-\Phi(\upsilon)|+|f_2-\upsilon|+|\Phi(\upsilon)-\upsilon|
\leq C(R)(1+C(R))\Upsilon^2(\rho)+|\Phi(\upsilon)-\upsilon|.
$$
Из \eqref{ekv:3.2} и \eqref{ekv:3.3} следует
$|\Phi_2(\upsilon)|+|\Phi_3(\upsilon)|<C(R)\Upsilon^2(\rho)$,
поэтому
\begin{multline*}
|\Phi(\upsilon)-\upsilon|\leq|\upsilon-\Phi_0-\Phi_1(\upsilon)|+C(R)\Upsilon^2(\rho)\leq\\
\leq\left|4\intl_0^x\left(\intl_0^t\intl_0^s\sigma(t)\sigma(s)\sigma(\tau)
\cos2\rho t\cos2\rho s\sin2\rho \tau \,d\tau\, ds\right)dt\right|+\\
+\left|\frac1{\rho}\intl_0^x\intl_0^t\sigma(t)\sigma^2(s)\cos2\rho
t(1-\cos2\rho s)\,ds\,dt\right|\leq C(R)\Upsilon^2(\rho).
\end{multline*}
Так как двойной интеграл после изменения порядка интегрирования
можно представить в виде произведения обычных интегралов, то оба
интегральных слагаемых оцениваются величиной
$C(R)\Upsilon^2(\rho)$. Это доказывает неравенство \eqref{lt1} и,
тем самым, утверждение леммы.
\end{proof}

\begin{Lemma}\label{lem3.4}
Пусть \(\sigma(x)\in W_2^{\theta}\), \(0\leqslant\theta\leqslant 1\),
\(\|\sigma\|_{\theta}\leqslant R\). Тогда преобразование Фурье этой
функции
\begin{equation}\label{ekv:Fou}
    F(\rho)=\int\limits_0^{\pi}\sigma(x) e^{i\rho x}\,dx
\end{equation}
в полосе \(|\Im\rho|\leqslant\nu\) допускает оценку \(|F(\rho)|<CR
|\rho|^{-\theta}\), где \(C\) зависит только от \(\nu\).
\end{Lemma}
\begin{proof}
Это утверждение, по-видимому, является фольклором для
специалистов, но мы затрудняемся указать точную ссылку. Поэтому
представим доказательство.

При \(0\leqslant\theta\leqslant 1\) функция \(\sigma\in W_2^{\theta}\)
представима рядом по косинусам
\begin{equation}\label{ekv:sigma}
    \sigma(x)=\sum\nolimits_{k=0}^{\infty}c_k \cos kx,\qquad
    \text{где}\quad\sum\nolimits_{k=0}^{\infty}k^{2\theta}\,|c_k|^2\leq
    \dfrac{2}{\pi}R^2.
\end{equation}
Проведём оценку \(F(\rho)\) вне кругов \(|\rho-k|<1/2\), \(k=0,
\pm 1,\pm 2,\ldots\) (внутри этих кругов она будет справедлива по
принципу максимума для аналитических функций). Подставим ряд для
\(\sigma\) в интеграл и проинтегрируем почленно. Полагая
\(c_{-k}=c_k\), получим
\[
    |F(\rho)|\leqslant C\sum\limits_{k=-\infty}^{\infty}
    \dfrac{|c_k|}{|k-\rho|}\leqslant C\left(
    \sum\limits_{k=-\infty}^{\infty}|k|^{2\theta}|c_k|^2\right)^{1/2}
    \left(\sum\limits_{k=-\infty}^{\infty}\dfrac{1}{|k|^{2\theta}
    |k-\rho|^2}\right)^{1/2}.
\]
Теперь утверждение леммы следует из тривиальных интегральных оценок
\begin{align*}
    \int\nolimits_1^{\mu -1}
    \dfrac{dt}{t^{2\theta}(\mu-t)^2}&<C\mu^{-2\theta},&
    \int\nolimits_{\mu+1}^{\infty}\dfrac{dt}{t^{2\theta}(t-\mu)^2}&<
    C\mu^{-2\theta}.
\end{align*}
Лемма доказана.
\end{proof}

\begin{Lemma}\label{lem3.5}
Пусть последовательность \(\{\rho_n\}_1^{\infty}\) такова, что
\(|\rho_n-n|<\delta<1/4\). Тогда при \(0\leq\theta<1/2\) оператор
\(T_x: W_2^{\theta}\to\ell_2^{\theta}\), определённый равенством
\begin{align*}
    T_xf&=\{c_n\}_1^{\infty},&c_n(x)&=\int\nolimits_0^x f(t)
    e^{i\rho_n t}\,dt,
\end{align*}
ограничен и его норма зависит только от \(\delta\) и \(\theta\).
\end{Lemma}
\begin{proof}
Согласно теореме Кадеца~\cite{Ka}, система \(\{e^{i\rho_n
x}\}_1^{\infty} \cup 1\cup\{e^{-i\rho_n x}\}_1^{\infty}\) образует
базис Рисса в пространстве \(L_2[0,2\pi]\), если
\(|\rho_n-n|<\delta<1/4\), при этом норма оператора, переводящего
эту систему в произвольную ортонормированную систему, зависит
только от \(\delta\). Следовательно, оператор
\(T_{2\pi}:L_2[0,2\pi]\to \ell_2\) ограничен и его норма
\(\leqslant C(\delta)\). Если \(f(x)\in\Wo_2^1[0,2\pi]\) (так
обозначаем подпространство функций в \(W_2^1[0,2\pi]\),
аннулирующихся на концах отрезка \([0,2\pi]\)), то интегральные
представления для чисел \(c_n(2\pi)\) можно проинтегрировать по
частям. Тогда с учетом оценки \(|\rho_n|>n/2\) получаем, что
оператор
\[
    T_{2\pi}:\Wo_2^1[0,2\pi]\to\ell_2^1
\]
также ограничен и его норма \(\leqslant 2C(\delta)\). Если
\(\Wo_2^{\theta}= [\Wo_2^{\theta},L_2]_{\theta}\) "---
промежуточное пространство, то, в силу теоремы об интерполяции
(см., например, \cite{Trl}), норма оператора
\(T_{2\pi}:\Wo_2^{\theta}\to\ell_2^{\theta}=[\ell_2^1,\ell_2]_{\theta}\)
не превосходит \(2^{\theta}C(\delta)\). Теперь заметим
(см.~\cite{BIN}, \cite{Trl}), что при \(0\leqslant\theta<1/2\)
пространства \(\Wo_2^{\theta}\) и \(W_2^{\theta}\) совпадают, а
оператор умножения \(\chi_x\) на характеристическую функцию
отрезка \([0,x]\) ограничен в пространстве \(W_2^{\theta}\). Из
равенства \(T_x=T_{2\pi}\chi_x\) получаем утверждение леммы.
\end{proof}

\begin{Lemma}\label{lem3.6}
Пусть $\sigma(x)\in L_2(0,\pi), $ $\|\sigma\| \leqslant  R. $
Тогда все корни функции $s(\pi,\rho) $ лежат в полосе $|\Im\rho|
<4e^{2\pi R}$.
\end{Lemma}
\begin{proof}
Известно \cite{HM4}, что справедливo представление
\begin{equation}\label{eq2.2}
\rho s(x,\rho)=\sin(\rho x)+\intl_0^x K(x,t)\sin (\rho t)dt,
\end{equation}
где ядро $K(x,t)\in L_2(0,\pi)$ при каждом фиксированном
$x\in[0,\pi]$, в частности, $\| K((\pi,x)\| \leq C(\sigma)$. В
работе \cite{HM4} явно не указывалась зависимость $C$  от
$\sigma$, но из анализа доказательств (см. также \cite{HM1} и
\cite{HM3}) легко видеть, что $C< e^{\pi R}.$ Тогда интегральное
слагаемое в \eqref{eq2.2} при $\rho =\mu+i\nu$ оценивается
величиной
$$
\int_0^{\pi}|K(\pi,t)|\, e^{\nu t}\, dt \leq e^{\pi
R}\left(\int_0^{\pi} e^{2\nu\pi}\, dt\right)^{1/2} <
e^{\pi(R+\nu)}\frac1{\sqrt{2\nu}}.
$$
Поскольку $|\sin\pi\rho|> (1/2)e^{\pi\nu}$,  то при $\nu=|\Im\rho|
>4e^{2\pi R}$ уравнение \eqref{eq2.2} корней не имеет.
\end{proof}

 Теперь завершим доказательство теоремы~\ref{thm3.1}. Из
леммы~\ref{lem3.2} следует, что корни функций \(s(\pi,\rho)\) и
\(\sin\theta(\pi,\rho)\) совпадают с числами \(\pm\rho_n\).
Оператор умножения \(\chi_x\) на характеристическую функцию
отрезка \([0,x]\) ограничен в пространствах \(W_2^{\theta}\) при
\(0\leqslant\theta< 1/2\) константой, зависящей от $\theta$, но не
зависящей от $x\in[0,\pi]$. Поэтому из леммы~\ref{lem3.4}
вытекает, что неравенство \eqref{ekv:ineq} выполнено для всех
точек $\rho$ полосы $|Im\,\rho|\leq\nu$ при
$$
|\rho|>\left[4\eps^{-1}RC(\nu)(1+64R^2\varkappa^2)^2+
R^2\varkappa+R^3\varkappa^2\right]^{1/\theta}=r_0=r_0(\nu,\theta,R).
$$
Это означает, что для всех таких точек $\rho$ справедливо
утверждение леммы \ref{lem3.3}: функция $\theta(\rho,x)$ корректно
определена и справедливы представления \eqref{ekv:pred},
\eqref{eq:pred1}. На границе круга $|\rho-n|\leq CRn^{-\theta}$,
согласно лемме \ref{lem3.4}, при всех $n>n_0=n_0(R)$ справедлива
оценка
$$
|f(\pi,\rho)|\leq C\Upsilon(\rho)<CR|\rho|^{-\theta}=|\pi\rho-\pi
n|.
$$
Из представления \eqref{ekv:pred} и теоремы Руше получаем, что
уравнение $\theta(\pi,\rho)=\pi n$ имеет в этом круге ровно один
корень. Как доказано в нашей работе \cite{SS2}( см. теорему 2.8),
этот корень будет иметь номер $n$. Выбором числа $n_0$ мы можем
добиться выполнения условия $CRn^{-\theta}<\delta<1/4$. Из
леммы~\ref{lem3.5} тогда имеем
$\{\Upsilon(\rho_n)\}_{n_0}^\infty\in l_2^\theta$ и
$\|\Upsilon(\rho_n)\|_\theta=O(1)$. Здесь и далее через \(O(1)\)
обозначается величина, которая оценивается константой, зависящей
от \(R\), \(\theta\), но не зависящей от \(\sigma\) в шаре
\(\|\sigma\|_{\theta}\leqslant R\). В силу леммы~\ref{lem3.6}
первые $n_0$ нулей лежат в прямоугольнике ширины $<C(R)$ и длины
$<n_0+1$. Поэтому такая же оценка сохранится для
последовательности $\{\Upsilon(\rho_n)\}_{1}^\infty$.

Положим теперь
$$
s_{2n}:=\rho_n-n=-\frac12b_{2n}+\al_{2n},\quad\text{где}\
b_k=\frac2\pi\intl_0^\pi\sigma(t)\sin ktdt.
$$
Мы уже доказали, что
 $|s_{2n}|=O(n^\theta)$. Уравнение
\(\theta(\pi,\rho_n)=\pi n\) перепишем в виде
\begin{multline}\label{dop}
    n\pi=\pi(n+(-1/2)b_{2n}+\alpha_{2n})+\intl_{0}^{\pi}\sigma(t)
    \sin 2(n+s_{2n})t\,dt+\\
    +2\intl_0^\pi\intl_0^t\sigma(t)\sigma(s)\cos2(n+s_{2n})t\sin2(n+s_{2n})s\,dsdt+r_n,
\end{multline}
где $r_n \leq \Upsilon^2(\rho_n) +O(n^{-1})$. Поэтому $\{r_n\}\in
l_2^{2\theta}$ и $\|\{r_n\}\|_{2\theta}=O(1)$ (заметим, что
последовательность $\{n^{-1}\}_1^{\infty}$ принадлежит
пространству \(\hat{l}_2^{2\theta}\) при
\(0\leqslant\theta<1/4\)). Используя тригонометрические формулы,
первый интеграл представим в виде
$$
    -\int\nolimits_0^{\pi}\sigma(t)\,(1-\cos 2s_{2n}t)
    \sin(2nt)\,dt+\pi/2\,b_{2n}+\int\nolimits_0^{\pi}\sigma(t)\,(\sin 2s_{2n} t)
    \cos 2nt\,dt
$$
Слагаемое $\pi/2\,b_{2n}$ сократится с таким же слагаемым в
\eqref{dop}, а для оценки оставшихся двух интегралов разложим
функции \(1-\cos 2s_{2n}t\) и \(\sin 2s_{2n}t\) в ряд по степеням
\(2s_{2n}t\). Получим, что эти два интеграла оцениваются величиной
\begin{gather}\label{ekv:3.5}
    |s_{2n}|^2\sum\limits_{k=1}^{\infty}
    \dfrac{1}{k!}|c_{k\,n}|+|s_{2n}|\sum\limits_{k=1}^{\infty}
    \dfrac{1}{k!}|b_{k\,n}|,\\
    \intertext{где}\notag
    \begin{aligned}
    c_{k\,n}&=\int\limits_{0}^{\pi}(2t)^k\,\sigma(t)\cos(2nt)\,dt,&
    b_{k\,n}&=\int\limits_{0}^{\pi}(2t)^k\,\sigma(t)\sin(2nt)\,dt.
    \end{aligned}
\end{gather}
Здесь, не ограничивая общности, мы считаем \(|s_{2n}|\leqslant
1\), так как это неравенство может нарушаться только для индексов,
число которых \(\leqslant C(R)\). Из определения нормы в
пространствах \(W_2^{\theta}\) при \(0<\theta<1/2\) следует
\begin{equation}\label{ekv:3.6}
    \sum\nolimits_{n=1}^{\infty} n^{2\theta}|c_{k\,n}|^2+
    \sum\nolimits_{n=1}^{\infty} n^{2\theta}|b_{k\,n}|^2=
    \pi\|(2t)^k\sigma\|^2_{\theta}<(2\pi)^{2k+1}\|\sigma\|^2_{\theta}.
\end{equation}
Из оценок~\eqref{ekv:3.5} и~\eqref{ekv:3.6} следует, что при
фиксированных \(k\) последовательности
\(\{|s_{2n}|b_{k\,n}\}_1^{\infty}\) и \(\{|s_{2n}|^2
c_{k\,n}\}_1^{\infty}\) принадлежат пространству
\(\ell_2^{2\theta}\) и их нормы ограничены числом
\(C(2\pi)^k/k!\).  Аналогичным приемом оценка двойного интеграла
из~\eqref{dop} сводится к оценке нормы $\|\{w_{n}\}\|_{2\theta}$,
где
$$
w_n=\intl_0^\pi\intl_0^t\sigma(t)\sigma(s)\cos2nt\sin2nsdsdt.
$$
Оценка $l_2^{2\theta}$ -- нормы этой последовательности проведена
в нашей работе \cite[предложение 3.8]{SS2}. Итак, числа
\(\{\alpha_{2n}\}\) в \eqref{dop} таковы,
что\(\{\alpha_{2n}\}_1^{\infty}\in l_2^{2\theta}\), и норма этой
последовательности есть $O(1)$.

Аналогично доказывается утверждение теоремы для последовательности
\(\{\alpha_{2n-1}\}_1^{\infty}\). Нужно только учесть, что \(\rho
s^{[1]}(x,\rho)=r(x,\rho)\cos\theta(x,\rho)\) (см.~лемму~2.5
работы~\cite{SS2}). Теорема доказана.
\end{proof}
\begin{Note}
В случае \(\sigma\in W^{\theta}\), \(0\leqslant\theta<1/2\) в
работе~\cite{SS2} доказано, что \(\{\alpha_k\}_1^{\infty}\in
 l_p^{2\theta}\) не только при \(p=2\), но при всех \(p>1\).
Совершенствуя метод статьи~\cite{SS2} можно получить б\'{о}льшее:
норма последовательности \(\{\alpha_k\}_1^{\infty}\) в \
 \(l_p^{2\theta}\), \(p>1\), \(0\leqslant\theta<1/2\) оценивается
константой \(C\), зависящей только от \(\theta\), \(p\), \(R\), но
не зависящей от \(\sigma\) в шаре \(\|\sigma\|_{\theta,\,
p}\leqslant R\), т.е. теорема \eqref{thm3.1} допускает
существенное усиление. Но это требует более длительной работы.
\end{Note}

\medskip

\section*{\bf \S3. Случай $\theta=1$.}\refstepcounter{section}

\medskip

При \(\theta=1\) имеем классический потенциал \(\sigma'(x)=q(x)\in
L_2\). Естественно, этот случай наиболее прост. Основная теорема в
этом случае эквивалентна следующему утверждению.

\begin{Theorem}\label{thm4.1}
При \(\sigma(x)\in W_2^1[0,\pi]\) справедливы формулы
\begin{align*}
    s_{2k}&=\dfrac{h_1}{2k}-\dfrac{a_{2k}}{2(2k)}+
    \dfrac{\alpha_{2k}}{k^2},&
    s_{2k-1}&=\dfrac{g_1}{2k-1}-\dfrac{a_{2k-1}}{2(2k-1)}+
    \dfrac{\alpha_{2k-1}}{k^2},
\end{align*}
где
\begin{align*}
    h_1&=\dfrac{\sigma(\pi)-\sigma(0)}{\pi}=\dfrac1\pi\int\limits_0^{\pi}
    q(x)\,dx,& g_1&=\dfrac{-\sigma(0)-\sigma(\pi)}\pi,&
    a_p&=\dfrac{2}{\pi}\int\limits_0^{\pi}q(t)\,\cos pt\,dt,&
    q(x)&=\sigma'(x),
\end{align*}
а числа \(\alpha_n\) таковы, что \(\sum_{n=1}^{\infty}|\alpha_n|^2
\leqslant C(R)\), где постоянная \(C(R)\) зависит только от
\(R=\|q\|= \|q\|_{L_2}\leq \|\sigma\|_1\).
\end{Theorem}
\begin{proof}
Эквивалентность этого утверждения и основной теоремы при \(\theta=1\)
становится ясной, если равенства~\eqref{T} проинтегрировать по частям.
Тогда
\begin{align*}
    b_{2k}&=\dfrac{h_1}{2k}-\dfrac{a_{2k}}{2(2k)},&
    b_{2k-1}&=\dfrac{g_1}{2k-1}-\dfrac{a_{2k-1}}{2(2k-1)}.
\end{align*}
Остаётся провести оценку нормы последовательности \(\{\alpha_k\}\)
в пространстве \(l_2\).

Как прежде, будем работать с числами $s_{2k}=\sqrt{\la_k}-k$.
Изменения, которые нужно сделать при работе с числами  $s_{2k-1}$,
очевидны.

Пусть $s(x,\rho)$ -- решение уравнения
$$
-y''+q(x)y=\rho^2y,\quad \sigma(x)=\intl_0^x q(t)dt,
$$
с начальными условиями $s(0,\rho)=0$, $s'_x(0,\rho)=\rho$ (здесь
удобнее изменить начальные условия для функции $s$, определенной
ранее в п.2). Используя метод последовательных приближений,
представим это решение в виде
\begin{equation}\label{S}
s(x,\rho)=\suml_{n=0}^\infty S_n(x,\rho),
\end{equation}
где
$$
S_0(x,\rho)=\sin(\rho x),\quad S_n(x,\rho)=
\intl_0^x\frac{\sin\rho(x-t)}{\rho}q(t)S_{n-1}(t,\rho)dt,\ n=1,2,
\dots.
$$
В силу оценки
\begin{equation}\label{n8}
|S_n(x,\rho)|\leq\frac{\|q\|^nx^{n/2}}{\sqrt{n!}|\rho|^n}e^{nx|Im\rho|},\quad
n=1,2,\dots
\end{equation}
(легко проверяемой по индукции) ряд сходится.  Как и ранее,
обозначения $O(\rho^{-n})$ или $O(k^{-n})$  используем для
выражений, абсолютная величина которых допускает оценки $\leq
C|\rho|^{-n}$ или $\leq C k^{-n},\ k\geq1$, где постоянная $C$
зависит только от $R=\|q\|$. В первом случае всегда предполагаем,
что оценка выполняется в полосе $|\Im\rho|\leq C=C(R)$.

Из оценки \eqref{n8} следует, что
\begin{equation}\label{n9}
\left|\suml_{n=3}^\infty
S_n(x,\rho)\right|=O(\rho^{-3})\quad\text{при}\ 0\leq x\leq\pi.
\end{equation}
Далее,
\begin{equation}\label{n10}
S_1(x,\rho)=-\frac{\cos(\rho x)}{2\rho}\intl_0^x
q(t)dt+\intl_0^x\frac{\cos\rho(x-2t)}{2\rho}q(t)dt,
\end{equation}
\begin{multline}\label{n11}
S_2(x,\rho)=\\ -\frac{\sin(\rho
x)}{(2\rho)^2}\intl_0^x\sigma(t)q(t)dt-
\intl_0^x\frac{\sin\rho(x-2t)}{(2\rho)^2}\sigma(t)q(t)dt+
\intl_0^xq(t)\sin\rho(x-t)\intl_0^t\frac{\cos\rho(t-2\xi)}{2\rho^2}q(\xi)d\xi
dt.
\end{multline}
Из выписанных соотношений следует, что
$$
s(\pi,\rho)=\sin\rho\pi+O(\rho^{-1}).
$$
Поэтому нули $\rho_k$ этой функции удовлетворяют соотношению
\begin{equation}\label{n12}
\rho_k=k+s_{2k},\qquad s_{2k}=O(k^{-1}).
\end{equation}
Заметим, что
$$
\intl_0^t\frac{\cos\rho(t-2\xi)}{2\rho^2}q(\xi)d\xi=
\frac{\cos\rho t}{2\rho^2}
\intl_0^\pi\cos(2\rho\xi)q(\xi)\chi_t(\xi)d\xi+ \frac{\sin\rho
t}{2\rho^2}\intl_0^\pi q(\xi)\chi_t(\xi)\sin(2\rho\xi) d\xi,
$$
где $\chi_t(\xi)$ -- характеристическая функция интервала $[0,t]$.
Разлагая функции $\cos2(k+s_{2k})\xi$ и $\sin2(k+s_{2k})\xi$ в
сумму произведений, получаем, что при $\rho=\rho_k=k+s_{2k}$ оба
слагаемых в правой части этого выражения имеют вид
\begin{equation}\label{n13} \{\al_k O(k^{-2})\}, \qquad
\text{где}\ \sum|\al_k|^2\leq\|q\chi_t\|^2\leq\|q\|^2.
\end{equation}
Следовательно, третье слагаемое в правой части \eqref{n11} при
$x=\pi$  и $\rho=\rho_k$ имеет такой же вид. Тогда
\begin{equation}\label{n14}
S_2(\pi,\rho_k) =\gamma_kk^{-2}, \qquad \sum |\gamma_k|^2 \leq
C(R).
\end{equation}
(здесь мы подразумеваем, что получение такого представления для
первых двух слагаемых в правой части \eqref{n11} проще, чем для
третьего).

Далее, с учетом введенных обозначений из \eqref{n10} имеем
\begin{equation}\label{S_1}
S_1(\pi,\rho_k)=(-1)^{k+1}\pi\left(\frac{h_1}{2k}-\frac{a_{2k}}{2(2k)}\right)
+\gamma'_kk^{-2},\quad\text{где}\ \sum|\gamma'_k|^2\leq C(R).
\end{equation}
 Запишем числа $\rho_k$ в виде
$$
\rho_k=k+\frac{h_1}{k}-\frac{a_{2k}}{2(2k)}+\frac{\delta_k}{k^2},
$$
где $\delta_k$ -- некоторые числа, которые в силу \eqref{n12}
подчинены оценке $\delta_k=O(k)$. Подставим это выражение в
равенство
\begin{equation}\label{S=0}
\sin\rho_k\pi+S_1(\pi,\rho_k)+S_2(\pi,\rho_k)+
O\left(k^{-3}\right)=0.
\end{equation}
(здесь мы учли, что остаток ряда \eqref{n9} есть $O(k^{-3}))$.
Заметим, что $\left(\delta_kk^{-2}\right)^3 = O(k^{-3})$, поэтому
$$
\sin\rho_k\pi
=(-1)^k\pi\sin\left(\frac{h_1}{2k}-\frac{a_{2k}}{2(2k)}+\frac{\delta_k}{k^2}\right)=
(-1)^k\pi\left(\frac{h_1}{2k}-\frac{a_{2k}}{2(2k)}+\frac{\delta_k}{k^2}\right)+
O\left(k^{-3}\right).
$$
Подставив выражения  \eqref{n14} и \eqref{S_1} в \eqref{S=0},
получим $(-1)^k\pi\delta_k=-\gamma_k-\gamma_k'+O(k^{-1})$.  Это
равенство доказывает основную теорему при $\theta=1$.
\end{proof}

\medskip

\section*{\bf \S4. Случай целых $\theta =m>1$.}\refstepcounter{section}

\medskip

Утверждение основной теоремы при целых $\theta=m>1$ получается из
следующего утверждения.
\begin{Theorem}\label{thm5.1}
Пусть $\sigma(x)\in W_2^m$, т.е. $q(x)\in W_2^{m-1}$, $m\geq1$.
Тогда для собственных значений $\la_k$ и $\mu_k$ операторов $L_D$
и $L_{DN}$ справедливы следующие формулы. При нечетном $m=2s+1$
\begin{gather}\label{n23}
\la_k^{1/2}=k+\dfrac{h_0}{(2k)}+\dfrac{h_1}{(2k)^3}+\dots+\dfrac{h_s}{(2k)^{2s+1}}-
(-1)^s\dfrac{a_{2k}}{2(2k)^{2s+1}}+\dfrac{\alpha_{2k}}{k^{2s+2}},\\
\notag
\mu_k^{1/2}=k-\dfrac12+\dfrac{g_0}{(2k-1)}+\dfrac{g_1}{(2k-1)^3}+\dots+
\dfrac{g_s}{(2k-1)^{2s+1}}-
(-1)^s\dfrac{a_{2k-1}}{2(2k-1)^{2s+1}}+\dfrac{\alpha_{2k}}{k^{2s+2}}.
\end{gather}
При четном $m=2s$
\begin{gather}\label{n24}
\la_k^{1/2}=k+\dfrac{h_0}{(2k)}+\dfrac{h_1}{(2k)^3}+\dots+\dfrac{h_s}{(2k)^{2s+1}}-
(-1)^s\dfrac{b_{2k}}{2(2k)^{2s}}+\dfrac{\alpha_{2k}}{k^{2s+1}},\\
\notag
\mu_k^{1/2}=k-\dfrac12+\dfrac{g_0}{(2k-1)}+\dfrac{g_1}{(2k-1)^3}+\dots+\dfrac{g_s}{(2k-1)^{2s+1}}-
(-1)^s\dfrac{b_{2k-1}}{2(2k-1)^{2s}}+\dfrac{\alpha_{2k-1}}{k^{2s+1}}.
\end{gather}
В этих формулах $ \{\alpha_l\}_1^\infty$ --- последовательность из
$l_2$ и ее $l_2$-норма оценивается постоянной $C$,  зависящей от
$R$, но не зависящей от $q$  в шаре $\|q\|_{m-1} \leq R$.Числа
$a_l=a_l(q^{(m-1)})$ и $b_l=b_l(q^{(m-1)})$ определяются формулами
$$
a_l=\frac2\pi\intl_0^\pi q^{(m-1)}(t)\cos ltdt,\quad
b_l=\frac2\pi\intl_0^\pi q^{(m-1)}(t)\sin ltdt, \qquad
l=1,2,\dots,
$$
а числа $h_j, g_j, \ 0\leq j\leq s$   являются непрерывными
ограниченными функционалами от $\sigma\in W_2^m$  и их линейные
части $h^0_j,g^0_j$ при $j\leq s-1$ выражаются формулами
$$
h_j^0=(-1)^{j+1}\pi^{-1}[\sigma^{(2j)}(\pi)-\sigma^{(2j)}(0)],\quad
g_j^0=(-1)^{j}\pi^{-1}[\sigma^{(2j)}(\pi)+\sigma^{(2j)}(0)].
$$
В случае $m=2s+1$ эти формулы сохраняются при $j=s$,  а в случае
$m=2s$ функционалы $h^0_s, g^0_s$ обращаются в ноль.
\end{Theorem}
\begin{proof}
Чтобы понять эквивалентность утверждений основной теоремы и
теоремы \ref{thm5.1} нужно равенства \eqref{T}  проинтегрировать
по частям $m$ раз.  Получим, что четные координаты
последовательности $T(\sigma)$ выражаются формулами \eqref{n23} в
которых функционалы $h_j$ заменяются на $h^0_j$, а остаточные
члены $\alpha_{2k}k^{-m-3}$ полагаются равными нулю. Аналогичное
утверждение справедливо для нечетных координат.   Но тогда
непрерывность функционалов $h_j-h^0_j$  и $g_j-g^0_j$  от $\sigma
\in W_2^m$  и равномерная оценка $l_2$-нормы последовательности
$\{\alpha_l\}$ влечет оценку \eqref{1.4}.

Для доказательства теоремы будем вновь использовать представление
\eqref{S}. Нам нужно иметь точную информацию только о слагаемых
$S_j$ при $j\leq m+1$, для суммы остальных слагаемых согласно
\eqref{n8} получаем оценку
\begin{equation} \label{rem}
\Big|\suml_{n=m+2}^{\infty}S_{n}(x,\rho)\Big|=O(\rho^{-m-2}).
\end{equation}

Введем обозначения
$$
\nu_{2s}(x,\rho):=(2\rho)^{-2s}\sin\rho x,\quad
\nu_{2s+1}(x,\rho):=(2\rho)^{-(2s+1)}\cos\rho x,\quad s=0,1,\dots,
$$
$$
(\pm)_j=\left\{\begin{array}{ll} -1,\quad\text{если}\ j=4s,\
4s+1,\\\phantom{-}1,\quad\text{если}\quad j=4s+2,\
4s+3.\end{array}\right.
$$
Основой дальнейшего доказательства является следующее предложение.

\smallskip
\begin{Lemma}\label{lem2}
Пусть $q(x)\in W_2^{m-1}, \sigma(x)=-\int_x^\pi q(t)\,dt, m\geq2$.
Тогда
\begin{equation}\label{n15}
S_1(x,\rho)=\suml_{j=1}^m\nu_j(x,\rho)f_{1,j}(x)+
(\pm)_{m+1}\intl_0^{x}\nu_{m}(x-2t,\rho)q^{(m-1)}(t)dt,
\end{equation}
\begin{equation}\label{n16}
S_2(x,\rho)=\suml_{j=1}^m\nu_j(x,\rho)f_{2,j}(x)+
(\pm)_{m+2}\intl_0^{x}\nu_{m+1}(x-2t,\rho)(q(t)\sigma(t))^{(m-1)}dt+O(\rho^{-m-2}),
\end{equation}
\begin{equation}\label{n17}
S_p(x,\rho)=\suml_{j=1}^m\nu_j(x,\rho)f_{p,j}(x)+O(\rho^{-m-2}),\quad
p=3,\dots,m+1.
\end{equation}
где
\begin{equation}\label{n18}
f_{1,j}(x)=(\pm)_j(\sigma^{(j-1)}(x)-(-1)^{j-1}\sigma^{(j-1)}(0)),
\end{equation}
а при $p\geq2$ функции $f_{p,j}(x)$ определяются рекуррентными
соотношениями

\begin{equation}\label{n19}
f_{p,j}(x)=(-1)^j\int_0^x q(t)f_{p-1,j-1}(t)\,
dt-\sum_{s=j}^m(\pm)_s(\pm)_j\left([qf_{p-1,s}]^{(j-s-2)}(x)-
[qf_{p-1,s}]^{(j-s-2)}(0)\right).
\end{equation}
При этом $f_{p,j}(x)\in W_2^{m+p-j-1}[0,\pi]$,  т.е. при каждом
фиксированном $x$ $f_{p,j}(x)$ являются непрерывными ограниченными
функционалами от $\sigma\in W^m_2$.
\end{Lemma}
\begin{proof}
Далее будем использовать тригонометрические тождества
\begin{equation}\label{n20}
\frac{\sin\rho(x-t)}{\rho}\nu_j(t,\rho)=(-1)^{j+1}\nu_{j+1}(x,\rho)+
\nu_{j+1}(x-2t,\rho)
\end{equation}
и тождества
\begin{multline}\label{n21}
\intl_0^x\nu_j(x-2t,\rho)f(t)dt=(-1)^j\intl_0^xf(t)d\nu_{j+1}(x-2t,\rho)=\\
=\nu_{j+1}(x,\rho)(f(x)-(-1)^jf(0))+(-1)^{j+1}\intl_0^x\nu_{j+1}(x-2t,\rho)f'(t)dt,
\end{multline}
справедливые для функций $f(t)\in W_2^1[0,\pi]$. В принятых
обозначениях имеем
$$
S_1(x,\rho)=\intl_0^x\frac{\sin\rho(x-t)}{\rho}\nu_0(t,\rho)q(t)dt=
-\nu_1(x,\rho)\left[\sigma(x)-\sigma(0)\right]+\intl_0^x\nu_1(x-2t,\rho)q(t)dt.
$$
Интегрируя по частям $m-1$ раз второе слагаемое в правой части
этого равенства и используя тождества \eqref{n21} получаем
представление \eqref{n15}, где функции $f_{1,j}$ определены
\eqref{n18}. Чтобы получить представления для других функций
$S_j$, воспользуемся тождеством
\begin{multline}\label{n22}
\intl_0^x\frac{\sin\rho(x-t)}{\rho}\nu_s(t,\rho)f(t)dt=\\
=(-1)^{s+1}\nu_{s+1}(x,\rho)\intl_0^xf(t)\, dt+
\intl_0^x\nu_{s+1}(x-2t,\rho)f(t)dt=\\
=(-1)^{s+1}\nu_{s+1}(x,\rho)\int_0^xf(t)\, dt-
\suml_{j=s+2}^{m+1}\nu_j(x,\rho)(\pm)_s(\pm)_j (f^{(j-s-2)}(x)-\\
-(-1)^{j-1}f^{(j-s-2)}(0))-
(\pm)_s(\pm)_{m+2}\intl_0^x\nu_{m+1}(x-2t,\rho)f^{(m-s)}(t)dt,
\end{multline}
которое получается из  соотношений \eqref{n20},
  если $f(x)\in W_2^{m-1}[0,\pi]$. Это тождество является основным
 для доказательства представлений \eqref{n16} и \eqref{n17}.
 Сначала заметим, что
 \begin{equation}\label{doubl}
 \intl_0^x\frac{\sin\rho(x-t)}{\rho}q(t)
 \intl_0^t\nu_m(t-2\xi,\rho)q^{(m-1)}(\xi)d\xi dt=O(\rho^{-m-2}).
 \end{equation}
 Это получается интегрированием по частям с учетом равенств
 $$
 \sin\rho(x-t)dt=\rho^{-1}d\cos\rho(x-t),\quad\nu_m(x,\rho)=O(\rho^{-m}).
 $$
 Теперь с учетом \eqref{n22}  и  \eqref{doubl} получаем
\begin{multline*}
 S_2(x,\rho)=\intl_0^x\frac{\sin\rho(x-t)}{\rho}q(t)
 \suml_{s=1}^m\nu_s(t,\rho)f_{1,s}(t)dt+O(\rho^{-m-2})=\\
 =-\suml_{s=1}^m\left((-1)^{s+1}\nu_{s+1}(x)\int_0^x q(t)f_{1,s}(t)\, dt+
 \suml_{j=s+2}^{m+1}\nu_j(t,\rho)(\pm)_s(\pm)_j([qf_{1,s}]^{(j-s-2)}(x)-\right.\\
 \left.-(-1)^{j-1}[qf_{1,s}]^{(j-s-2)}(0))-
 (\pm)_s(\pm)_{m+2}\intl_0^x\nu_{m+1}(x-2t,\rho)[q(t)\sigma(t)]^{(m-1)}dt\right)
 +O(\rho^{-m-2}).
\end{multline*}
 Здесь мы учли, что $f_{1,j}\in W_2^{m+1-j}[0,\pi]$ при $j\geq2$,
 поэтому при таких значениях индекса $j$ интегрирование по частям
 можно проводить до тех пор, пока интегральное слагаемое не
 станет равным $O(\rho^{-m-2})$. Меняя в последней формуле порядок
 суммирования и используя обозначения \eqref{n19}, получаем
 представление \eqref{n16}.

 Доказательство равенства \eqref{n17} проводится точно также, даже
 проще, так как не возникает проблемы с оценкой двойного интеграла \eqref{doubl}.
 Последнее утверждение леммы о принадлежности функций
 $f_{p,j}$ пространствам $ W_2^{m+p-j-1}$ является очевидным. Лемма доказана.
\end{proof}

Доказательство теоремы теперь завершим по индукции. Утверждение
уже доказано для $m=1$. Предположим, что формула \eqref{n23} верна
для $m-1=2s-1, s\geq 1$. Докажем, что формула  \eqref{n24}  верна
для $m=2s$.   Далее удобно использовать обозначение
$$
\Sigma_r = c_1k^{-1}+ c_2k^{-2}+\dots + c_rk^{-r},
$$
где $c_j$ --- некоторые числа. Через $\{\gamma_k\}$ обозначаем
различные последовательности, $l_2$-норма которых есть $O(1)$. С
учетом этих обозначений  формулу \eqref{n23} для $m-1$ запишем в
виде
$$
\lambda_k^{1/2}: =\rho_k
=k+\Sigma_{2s-1}-(-1)^ka_{2k}2^{-1}(2k)^{-2s+1}
+\gamma_{2k}k^{-2s} = :k+\Sigma_{2s-1}+\delta_kk^{-2s}.
$$
Здесь $a_{2k}=a_{2k}(q^{(m-2)}) =-(2k)^{-1}b_{2k}(q^{(m-1)})$, что
получается интегрированием по частям. Поэтому $l_2$-норма
неизвестной последовательности $\{\delta_k\}$ есть $O(1)$. Имея
такую априорную информацию, получаем
$$
\sin\rho_k\pi=(-1)^k\pi(\Sigma_{2s+1}
+\delta_kk^{-2s})+\gamma_kk^{-2s-1},\quad
\cos\rho_k\pi=(-1)^k(1+\Sigma_{2s+1}+\gamma_kk^{-2s-1}).
$$
Воспользовавшись равенством $(2\rho)^{2s}\nu_m(x-2t)=\sin\rho x\,
\cos2\rho t - \cos\rho x\, \sin 2\rho t$, оценкой \eqref{rem} и
представлениями \eqref{n15} - \eqref{n17}, находим
\begin{multline*}
0=s(\rho_k,\pi)=\Sigma_{2s+1}+(-1)^k\pi\delta_k
k^{-2s}+\\
+(\pm)_{m+1}(-1)^{k+1}\pi2^{-2s-1}k^{-2s}b_{2k}(q^{(m-1)})+
\gamma_kk^{-2s-1}+O(k^{-2s-2}),
\end{multline*}
Из этого равенства следует
$$
\delta_k=(\pm)_{m+1}b_{2k}(q^{(m-1)})2^{-2s-1} +\gamma_kk^{-1}
$$
что влечет формулу \eqref{n24}.Формально нужно провести еще один
шаг индукции: переход от четного $m=2s$ к нечетному $m=2s+1$. Но
это делается точно также. Теорема \ref{thm5.1} доказана.
\end{proof}
Тем самым утверждение основной теоремы  доказано при любом целом
$m\geq 1$.
\begin{Note}\label{noteA}
Из доказательства теоремы следует, что числа $h_1,\dots,h_s$
определяются рекуррентно и являются линейными комбинациями
коэффициентов разложений функций $S_j(\rho,\pi)$, $j=1,\dots,m+1$
по степеням $\rho^{-1}$, т.е. являются линейными комбинациями
значений функций $f_{j,\,p}(x)$ в точках $x=0$ и $x=\pi$. Это
замечание позволяет сделать вывод: функционалы $h_j, g_j$ являются
дифференцируемыми отображениями из пространства $W^{\theta}_2$ в
$\C$,  поскольку этим свойством обладают функционалы $f_{j,p}(0)$
и $f_{j,p}(\pi)$.
\end{Note}

\medskip

\section*{\bf \S5. Дифференцируемость отображения
$\Phi(\sigma)$.}\refstepcounter{section}

\medskip

Наблюдение о том, что собственные значения оператора
Штурма-Лиувилля аналитически зависят от вещественного потенциала,
принадлежит Боргу \cite{Bo}. Для вещественных потенциалов $q$ из
пространства $L_2$ эта идея получила существенное развитие в
работах Пошеля и Трубовица \cite{PT}.  Здесь сначала заметим, что
аналитичность сохраняется для потенциалов распределений из
пространств $W^{\theta}_2$ при $\theta \geq -1$. Кроме того,
вычисления производных мы будем проводить для комплексных
потенциалов.

Мы предполагаем, что читатель знаком с определениями {\it
производных по Фреше и Гато} для отображений $F:\,U\to H$, где $U$
-- открытое множество в $E$, а $E$ и $H$ -- банаховы пространства.
Далее используются факты, связанные с этими понятиями, которые
можно найти в книгах \cite{Di} и \cite{PT}. Линейный оператор из
$E$ в $H$, совпадающий с производной Фреше отображения $F$ в точке
$x$ будем обозначать через $d_xF$. В случае комплексных банаховых
пространств отображение $F:\,U\to H$, дифференцируемое по Фреше в
каждой точке $x\in U$ называется {\it аналитическим} в $U$.
Естественным образом определяется понятие {\it вещественного
аналитического отображения,} см. \cite{PT}. Напомним, что
отображение $F:\,U\to H$ называется {\it слабо аналитическим},
если для любых элементов $h\in E$, $x\in U$ и любого функционала
$L\in H^*$ скалярная функция $L(F(x+zh))$ дифференцируема в
комплексном смысле в некоторой малой окрестности нуля (зависящей
от элементов $x$, $h$). В случае гильбертова пространства $H$
слабая дифференцируемость эквивалентна дифференцируемости по Гато
координатных функций $(F(x),e_n)$, где $\{e_n\}$ --
ортонормированный базис в $H$. Мы будем использовать следующие
известные результаты, см.\cite {PT}.
\begin{Statement}\label{st:6.1}
Пусть $F:\,U\to H$ -- слабо аналитическое отображение множества
$U\in E$. Если $F$ локально ограничено в $U$, т.е. ограничено в
некоторой окрестности каждой точки $x\in U$, то $F$ --
аналитическое отображение в $U$.
\end{Statement}
\begin{Statement}\label{st:6.2}
Пусть $F:\,U\to H$ -- аналитическое отображение. Тогда справедлива
формула Коши
$$
F(x+zh)=\frac1{2\pi i}\intl_{|\zeta|=\eps}\frac{F(x+\zeta
h)}{\zeta-z}d\zeta,\quad |z|<\eps.
$$
Здесь $x\in U$, $h\in E$, а $\eps$ столь мало, что $x+zh\in U$ для
всех $|z|\leq\eps$. В частности,
$$
d_xF=\frac1{2\pi i}\intl_{|\zeta|=\eps}\frac{F(x+\zeta
h)}{\zeta^2}d\zeta,
$$
где $h$ -- произвольный элемент, такой, что $x+\zeta h\in U$ при
всех $|\zeta|\leq\eps$.
\end{Statement}
Найдем производные по Гато отображений $\la_k:\,L_2\to\C$, где
$\la_k=\la_k(\sigma)$ -- собственные значения оператора $L_D$. Для
собственных значений $\mu_k(\sigma)$ оператора $L_{DN}$ результаты
сохраняются.
\begin{Lemma}\label{lem:6.3}
Пусть $\la_k=\la_k(\sigma)$ -- простое собственное значение
оператора $L_D$, которому отвечает собственная функция
$y_k(x)=y_k(x,\sigma)$. Если $(y_k,\overline{y_k})=\int_0^\pi
y_k^2(x)dx\ne0$, то функция $la_k(\sigma):\,L_2\to\C$
дифференцируема по Фреше в окрестности точки $\sigma$ и ее
производная в этой точке равна
$$
[d_\sigma\la_k(\sigma)]h=-\frac{2(y_k'(x)y_k(x),\overline{h(x)})}{(y_k^2(x),1)}.
$$
Если $s_{2k}(\sigma): =\lambda^{1/2}_k -k,$ то
$$
[d_{\sigma}s_{2k}(\sigma)]h=
(1/2)\lambda^{-1/2}[d_{\sigma}\lambda_k(\sigma)]h.
$$
\end{Lemma}
\begin{proof}
Вспомним (\cite{SS1}), что
$$
\|y_k(x,\sigma)-y_k(x,\widetilde{\sigma})\|_1\to 0\quad\text{при}\
\|\sigma-\widetilde{\sigma}\|_{L_2}\to 0.
$$
Такое же соотношение выполняется для функций
$y_k^{[1]}(x,\sigma)$. Следовательно, функция $y_k'(x,\sigma)$
непрерывно зависит в $L_2$-норме от функции $\sigma\in L_2$.
Поэтому, если мы докажем существование производной по Гато в точке
$\sigma$ по направлению $h$:
$$
\lim\limits_{t\to0}\frac{\la_k(\sigma+th)-\la_k(\sigma)}t=-
\frac{2(y_k'y_k,\overline{h})}{(y_k^2,1)},
$$
то в силу непрерывности этой функции в малой окрестности точки
$\sigma$, она будет производной по Фреше.

Пусть $h\in L_2$ и $(\la_k(t),y_k(x,t))$ -- собственная пара
оператора $L_D$, отвечающая $\sigma+th$. Для краткости далее
опускаем индекс $k$. Запишем равенства
\begin{gather*}
-(y'(x,t)-(\sigma+th)y(x,t))'-(\sigma+th)y'(x,t)=
\lambda(t) y(x,t),\\
-(y'(x,0)-\sigma y(x,0))'-\sigma y'(x,0)= \lambda(0)y(x,0).
\end{gather*}
Умножим первое и второе равенства на \(y(x,0)\) и \(y(x,t)\),
соответственно, и выпишем разность. Проинтегрируем полученное
выражение (возьмем скалярное произведение с единичной функцией).
Проинтегрировав, где надо, по частям и воспользовавшись тем, что
обынтегрированные члены исчезают, получим
\begin{multline}\notag
  \big([(th(x)y(x,t))'-th(x)y(x,t)]y(x,0),1\big)
    =\\
    =(\lambda(t)-\lambda(0))\,\left(\left[y^2(x,0)+
    (y(x,t)-y(x,0))y(x,0)\right], 1\right).
\end{multline}
Заметим, что
$\left((h(x)y(x,t))'y(x,0),1\right)=-(h(x)y(x,t)y'(x,0),1)$.
Разделив полученное равенство на \(t\) и устремив \(t\) к нулю,
получим утверждение леммы. Для собственных значений оператора
$L_{DN}$ утверждение получается также, но надо более аккуратно
проводить интегрирование по частям.
\end{proof}
\begin{Lemma}\label{lem6.2}
Пусть $\sigma\in W_2^\theta$, $0<\theta<1/2$,
$\|\sigma\|_\theta\leq R$. Пусть $y_k=y_k(x,\sigma)$ --
собственные функции оператора $L_D$. Тогда найдется целое число
$N=N(R,\theta)$ и постоянная $C=C(R,\theta)$ такие, что
собственные значения $\la_k$ при $k\geq N$ простые и
\begin{equation}\label{N+1}
|y_k(x)-\sin kx|<Ck^{-\theta}.
\end{equation}
Такая же оценка сохраняется для собственных функций оператора
$L_{DN}$ с заменой $\sin kx$ на $\sin(k-1/2)x$.
\end{Lemma}
\begin{proof}
Это утверждение доказано в \cite[теорема 3.13]{SS2}, но с числами
$N$ и $C$, зависящими от $\sigma$ и $\theta$. Из результатов,
полученных в п.2 следует, что $N$ и $C$ зависят только от $R$ и
$\theta$.
\end{proof}
\begin{Theorem}\label{thm:6.5}
Пусть $\sigma\in W_2^\theta$, $\theta>0$ и $\sigma$ --
вещественная функция. Тогда
$F(\sigma):\,W_2^\theta\to\hat{l}_2^{\,\theta}$ -- вещественное
аналитическое отображение в окрестности точки $\sigma$.
\end{Theorem}
\begin{proof}
Для вещественной функции $\sigma(x)$ собственные функции $y_k(x)$
также вещественные, а потому $(y_k^2(x),1)>0$. В силу непрерывной
зависимости $y_k$ от $\sigma$ это неравенство сохраняется при всех
$k=1,\dots,\,N$ равномерно в некоторой окрестности точки $\sigma$.
Из оценки \eqref{N+1} получаем, что неравенства $(y_k^2,1)>0$
остаются в силе при всех $k\geq N+1$. Следовательно, координатные
функции пространства $l_2^{\,\theta}$ (в естественном базисе)
аналитичны в окрестности $\sigma$. Дифференцируемость координат
конечномерного пространства $\hat{l}_2^{\,\theta}\ominus
l_2^{\,\theta}$ следует из  замечания \ref{noteA}.
\end{proof}
Из аналитичности отображения $F(\sigma)$ следует аналитичность
отображения $\Phi(\sigma):\,W_2^\theta\to\hat{l}_2^{\,\tau}$, где
число $\tau$ определено в формулировке  основной теоремы. Конечно,
для комплексных $\sigma(x)$ дифференцируемости может не быть.

Обозначим через $P_N$ ортопроектор в $\hat{l}_2^{\,\theta}$,
аннулирующий координаты конечномерного пространства
$\hat{l}_2^{\,\theta}\ominus l_2^{\,\theta}$ и первые $N$
координат пространства $l_2^{\,\theta}$.
\begin{Theorem}\label{thm:6.6}
Пусть $0<\theta<1/4$. Тогда найдется число $N=N(R,\theta)$, такое,
что отображение
$F_N(\sigma)=P_NF(\sigma):\,W_2^\theta\to\hat{l}_2^{\,\theta}$
аналитично в шаре $\|\sigma\|\leq R$.
\end{Theorem}
\begin{proof}
Из неравенств \eqref{N+1} и леммы \ref{lem:6.3} следует, что при
$N\geq N(R,\theta)$ координаты отображения $F_N(\sigma)$ --
аналитические функции в шаре $\|\sigma\|\leq R$. Тогда в силу
предложения \ref{st:6.1} и теоремы \ref{thm3.1} отображение
$F_N(\sigma)$ аналитично в этом шаре.
\end{proof}
\begin{Theorem}\label{thm:6.7}
В условиях теоремы \ref{thm:6.6} отображение
$\Phi_N(\sigma)=P_N\Phi(\sigma):\,W_2^\theta\to\hat{l}_2^{\,2\theta}$
аналитическое и норма производной
$d_\sigma\Phi_N:\,W_2^\theta\to\hat{l}_2^{\,2\theta}$ ограничена
постоянной $C=C(R,\theta)$ для всех $\sigma$ из шара
$\|\sigma\|_\theta\leq R-1$. Это утверждение сохраняет силу при
всех целых $\theta=1,\,2,\dots$ (с заменой индекса $2\theta$ на
$\theta+1$).
\end{Theorem}
\begin{proof}
Из теорем \ref{thm3.1} и \ref{thm4.1} получаем ограниченность
отображения $\Phi_N:\,W_2^\theta\to\hat{l}_2^{\,2\theta}$ при
указанных значениях $\theta$. в шаре $\|\sigma\|_\theta\leq R$.
Так как $\Phi_N$ аналитическое, то из предложения \ref{st:6.2}
следует, что его производная ограничена в любом шаре радиуса $\leq
R-\eps$.
\end{proof}

\medskip

\section*{\bf \S6. Окончание доказательства основной теоремы}\refstepcounter{section}

\medskip

Воспользуемся следующим результатом, который сформулируем в
удобном для нас виде.
\begin{Statement}\label{st7.1}
Пусть $(E_1,E_0)$, $(H_1,H_0)$ пары банаховых пространств с
непрерывным вложением $E_1\subset E_0$, $H_1\subset H_0$. Пусть
$T$ -- нелинейное отображение из $E_0$ в $H_0$, отображающее $E_1$
в $H_1$ и удовлетворяющее следующему условию: существует
положительные возрастающие функции $C_0(R)$ и $C_1(R)$ такие, что
\begin{equation}\label{7.1}
\|\Phi\sigma-\Phi\widetilde{\sigma}\|_{H_0}\leq
C_0(R)\|\sigma-\widetilde{\sigma}\|_{E_0},
\end{equation}
если $\max\{\|\sigma\|_{E_0},\ \|\widetilde{\sigma}\|_{E_0}\}\leq
R$,
\begin{equation}\label{7.2}
\|\Phi\sigma\|_{H_1}\leq C_1(R)\|\sigma\|_{E_1}.
\end{equation}
Тогда $T$ отображает $E_\tau=[E_0,\,E_1]_\tau$ в
$H_\tau=[H_0,\,H_1]_\tau$ при всех $0\leq\tau\leq1$, причем
\begin{equation}\label{7.3}
\|\Phi\sigma\|_{H_\tau}\leq C_\tau(R)\|\sigma\|_{E_\tau},
\end{equation}
где $C_\tau(R)$ -- некоторая возрастающая функция от $R$.
\end{Statement}
\begin{proof}
 Доказательство этого предложения принадлежит Тартару \cite{Ta}.
\end{proof}
Воспользуемся этим предложением, положив $E_0=W_2^\theta$,
$H_0=\hat{l}_2^{\,\theta}$, где $\theta$ -- фиксированное число
$0<\theta<1/4$; $E_1=W_2^1$, $H_1=\hat{l}_2^{\,2}$, а в качестве
$\Phi$ возьмем отображение $\Phi_N$. Тогда из теоремы
\ref{thm:6.7} следует оценка \eqref{7.1}, а оценка \eqref{7.2}
была доказана в теореме \ref{thm4.1}. Следовательно, справедлива
оценка \eqref{7.3}, где $E_r=[W_2^\theta,\,W_2^1]_r$ и
$H_r=[\hat{l}_2^{\,2\theta},\,\hat{l}_2^{\,2}]_r$. Из теоремы о
реитерации (см., например \cite{Trl}) и теоремы об интерполяции
пространств $\hat{l}_2^{\,\theta}$ (см. \cite{SS3}) имеем
$E_r=W_2^{\theta+r(1-\theta)}$,
$H_r=\hat{l}_2^{\,2(\theta+r(1-\theta))}$. Следовательно,
отображение $\Phi_N:\,W_2^\theta\to\hat{l}_2^{\,2\theta}$ при
$0<\theta\leq1$ ограничено в любом шаре. Ограниченность
конечномерного отображения
$\Phi-\Phi_N:\,W_2^\theta\to\hat{l}_2^{\,2\theta}$ очевидна,
откуда получаем утверждение основной теоремы при $0<\theta\leq1$.

При $\theta\geq1$ нужно вновь воспользоваться предложением
\ref{st7.1}, но интерполировать нужно между пространствами $W_2^1$
и $W_2^m$, где $\theta\in[1,\,m]$. Этим заканчивается
доказательство основной теоремы.

Отметим, что доказав основную теорему, мы можем усилить результаты
п.5 об аналитичности отображений. А именно, имея локальную
ограниченность отображений $F$ и $\Phi$ при всех $\theta>0$, мы
можем усилить результаты п.5
\begin{Theorem}
Утверждение теорем \ref{thm:6.5} -- \ref{thm:6.7} сохраняют силу
при всех $\theta>0$. Отображение
$F(\sigma):\,W_2^\theta\to\hat{l}_2^{\,\theta}$ дифференцируемо в
некоторой комплексной окрестности вещественного шара
$\|\sigma\|_\theta\leq R$ при всех $\theta>0$ (а потому является
вещественно аналитическим отображением). Производная $F(\sigma)$ в
точке $\sigma=0$ равна $-(1/2)\, T $.
\end{Theorem}
\begin{proof}
Мы уже пояснили справедливость первых утверждений. Существование
производной по Фреше отображения $F(\sigma)$ в точке $\sigma=0$
при $\theta=0$ доказана в работе \cite{SS3}, причем показано, что
она равна $-(1/2)\,T$. Так как производная в точке $\sigma=0$
существует при всех $\theta\geq0$, пространства $W_2^\theta$
непрерывно вложены в $L_2$, то эта производная обязана совпадать с
$-(1/2)\, T$ при всех $\theta>0$.
\end{proof}

\end{document}